\documentclass[a4paper, 12pt, twoside, notitlepage]{amsart}

\usepackage{amsmath,amscd}
\usepackage{amssymb}
\usepackage{amsthm}
\usepackage{comment}
\usepackage{graphicx, xcolor}
\usepackage{mathrsfs}
\usepackage[ocgcolorlinks, linkcolor=blue]{hyperref}

\usepackage{bm}
\usepackage{bbm}
\usepackage{url}

\usepackage[utf8]{inputenc}
\usepackage{mathtools,amssymb}
\usepackage{esint}
\usepackage{dsfont}
\usepackage{relsize}
\usepackage{url}
\usepackage[shortlabels]{enumitem}
\usepackage{lineno}
\usepackage{enumitem}

\usepackage{tikz}
\usepackage{tikz-3dplot}
\usetikzlibrary{arrows.meta, calc}

\usepackage{subfigure}
\usepackage{array}
 \newcolumntype{C}[1]{>{\centering\arraybackslash}p{#1}}

 \usetikzlibrary{arrows,automata}

\usepackage{verbatim}
\usepackage{dsfont}
\numberwithin{equation}{section}

\allowdisplaybreaks

\graphicspath{{images/}}

\title[Survey on Broken Ray Transforms]{Survey on Broken Ray Transforms}
\author[S.R. Jathar]{Shubham R. Jathar}
\address{ Computational Engineering, School of Engineering Sciences,
Lappeenranta--Lahti University of Technology LUT, Lappeenranta, Finland}
\email {Shubham.Jathar@lut.fi}

\author[J. Railo]{Jesse Railo}
\address{
Computational Engineering, School of Engineering Sciences,
Lappeenranta--Lahti University of Technology LUT, Lappeenranta, Finland}
\email{Jesse.Railo@lut.fi}

\begin{document}

\begin{abstract}
We survey recent developments in the theory and applications of the broken ray transforms. Furthermore, we discuss some open problems.
\medskip


\end{abstract}
	\maketitle

\section{Introduction}

Let $(M, g)$ be a Riemannian manifold of dimension $n \geq 2$ with or without boundary. We consider a family of broken rays on $M$, which are piecewise smooth continuous curves that can either be geodesics or curves determined by some other dynamical equations. We denote the collection of such curves by $\Gamma$. For each $\gamma \in \Gamma$, the broken ray transform of a function $f \in C(M)$ is defined as
$$
R f(\gamma)=\int_\gamma f(\gamma).
$$
 One of the primary objectives is to analyze what information the line integrals of a function $f$ along the broken rays encodes about \(f\). The mathematical analysis of the broken ray transform involves investigating key properties such as injectivity, stability, reconstruction, partial data problems, and range characterization. A notable application is the partial data problem for the Calderón problem in electrical impedance tomography, where the invertibility of the broken ray transform plays a crucial role in reconstructing the conductivity of a medium from partial boundary measurements \cite{Kenig:Salo:2013}.

The points where the broken ray fails to be smooth are often referred to as the \textbf{points of reflection}. Depending on the location of these reflection points, broken rays can be classified into two primary categories.
\begin{itemize}
    \item \textbf{Boundary Reflections:} Reflection points \(z_i\) that lie on the boundary \(\partial M\). These represent cases where the ray interacts with the boundary and reflects according to the law of reflection, with the angle of incidence equal to the angle of reflection.
    \item \textbf{Interior Reflections (or Vertex Points):} Reflection points \(z_i\) that lie inside the manifold \(M\). These interior points represent situations where rays may change direction or ``break" at a vertex point within the domain, leading to transforms such as the V-line or star transforms. One physical motivation could be the scattering of particles.
\end{itemize}

\section{Boundary Reflections}
In this section, we consider the case where broken rays undergo reflections at the boundary of a compact manifold \((M, g)\) of dimension \(n \geq 2\). We assume that the boundary \(\partial M\) is decomposed into two distinct parts \(\mathcal{E}\) and \(\mathcal{R}\). A curve $\gamma$ on $M$ is called a broken ray if it starts and ends on $\mathcal E$, reflects on $\mathcal{R}$ according to the law of reflection (the angle of incidence equals the angle of reflection), and follows between reflections a geodesic (or other dynamical) trajectory in the interior of $M$.

\,

\noindent\textbf{Billiard-type Broken Rays.}
When both \(\mathcal{E}\) and \(\mathcal{R}\) are subsets of a connected boundary $\partial M$, the broken rays resemble billiard trajectories, where the curve reflects multiple times within the domain before returning to \(\mathcal{E}\) (see Figure \subref{fig:b},\subref{fig:c}). This case has been extensively studied quite extensively for some symmetric shapes like the disk and square. For example, in the case of a disk \(\mathbb{D} \subset \mathbb{R}^2\) with a non-empty open subset \(\mathcal{E}\) on its boundary, it has been established that if a function \(f \in C^{\alpha}(\overline{\mathbb{D}}; \mathbb{R})\) is quasianalytic in the angular variable, then \(f\) can be uniquely recovered from its broken ray transform \cite[Theorem 3]{Ilmavirta:2013}. On the other hand, there exist a smooth counterexample when data is over all closed billiards of the unit disk \cite[Proposition 33]{Ilmavirta:2015} 

In more general convex domains \(\Omega \subset \mathbb{R}^n\) with \(n \geq 2\), support theorems for the Radon line transform have been applied to establish local uniqueness \cite{Ilmavirta:2013}. Furthermore, for Riemannian manifolds with corners, the broken ray transform problem has been reduced to a geodesic ray transform on a related manifold through reflection, as demonstrated in \cite{Ilmavirta:2015}. These reductions employ established techniques from the theory of geodesic ray transforms.

Variants of this transform, such as the attenuated X-ray transform in a Euclidean square, have also been analyzed, leading to the results on injectivity and stability using microlocal analysis \cite{Hubenthal:2014} (see Figure \subref{fig:d}). These results were subsequently extended to higher dimensions \cite{Hubenthal:2015}. Additionally, the broken line transform in \(\mathbb{R}^2\) has been analyzed in \cite{Zhang:2020}, where it was identified as a Fourier integral operator, and conditions were determined under which certain singularities of a function are not recoverable due to the presence of conjugate points.
\begin{figure}[htbp]

\begin{tabular}{C{.3\textwidth}C{.3\textwidth}C{.3\textwidth}}
\subfigure [
] {
    \resizebox{0.3\textwidth}{!}{%
    \begin{tikzpicture}[scale=0.4, baseline={(current bounding box.south)}] 
  \draw[thick] (0,0) circle (3);

\draw[fill=red, opacity=0.6, thick] (0.7,0) ellipse (1.5 and 0.7);

\node at (3.3,-1) { $\mathcal{E}$};

\node at (0.7,0.2) { $\mathcal{R}$};

\draw[thick, ->] (-2, 2.23) -- (0.8, 0.7); 
\draw[thick, ->] (0.8, 0.7) -- (2.5, 1.8); 

\draw[thick, ->] (-3, -0.5) -- (2.6, -1.5);
    \end{tikzpicture}
    }\label{fig:a}
} & 
\subfigure [] {
  \resizebox{0.3\textwidth}{!}{%
   \begin{tikzpicture}[scale=0.4,baseline={(current bounding box.south)}]
  \draw[thick, blue] (180:3) arc[start angle=180, end angle=360, radius=3]; 
\draw[thick, red] (0:3) arc[start angle=0, end angle=180, radius=3]; 

\node at (2.8, -2) { $\mathcal{E}$}; 
\node at (3.2, 1.5) { $\mathcal{R}$};  

\draw[thick, ->] (-2.2,-2) -- (-1.2,2.75); 
\draw[thick, ->] (-1.2,2.75) -- (2.2,2); 

\draw[thick, ->] (2.2,2) -- (1.35, -2.7); 
    \end{tikzpicture}
    }\label{fig:b}
} 
 & 
\subfigure [] {
  \resizebox{0.3\textwidth}{!}{%
  \begin{tikzpicture}[scale=0.4, baseline={(current bounding box.south)}]
 \draw[thick] (0,0) circle (3);
\coordinate (A) at (0, 3);    
\coordinate (B) at (2.6, 1.5);  
\coordinate (C) at (2.6, -1.5); 
\coordinate (D) at (0, -3);    
\coordinate (E) at (-2.6, -1.5); 
\coordinate (F) at (-2.6, 1.5);  

\draw[->, thick] (A) -- (B);
\draw[->, thick] (B) -- (C);
\draw[->, thick] (C) -- (D);
\draw[->, thick] (D) -- (E);
\draw[->, thick] (E) -- (F);
\draw[->, thick] (F) -- (A);
    \end{tikzpicture}
    }\label{fig:c}
} 
\\

\subfigure [] {
  \resizebox{0.3\textwidth}{!}{%
    \begin{tikzpicture}[scale=0.6, baseline={(current bounding box.south)}]
\draw[thick] (0,0) rectangle (5,5);
\draw[thick, red] (0,0) -- (0,5);
\node at (-0.5, 2.5) { $\mathcal E$};
\filldraw[black] (0, 2) circle (2pt) node[left] { $x$};
\draw[thick, ->] (0, 2) -- (3, 5);  
\draw[thick, ->] (3, 5) -- (5, 2.5);  
\draw[thick, ->] (5, 2.5) -- (2.5, 0); 
\draw[thick, ->] (2.5, 0) -- (0, 3);     
\node at (2, 3.5) {$v$};
    \end{tikzpicture}
   }\label{fig:d}
} & 

\subfigure [] {
  \resizebox{0.3\textwidth}{!}{%
     \begin{tikzpicture}[scale=0.6, baseline={(current bounding box.south)}]

\draw[thick] (0,0) circle (3);

\begin{scope}[rotate=20] 
   
    \draw[thick, fill= red] 
        (-1.5, 0.6) .. controls (-1.2, 0.8) and (-0.8, 1.2) .. (0, 1.2) 
        .. controls (0.8, 1.2) and (1.2, 0.8) .. (1.5, 0.6)
        .. controls (1.6, 0.5) and (1.6, 0.4) .. (1.5, 0.3)
        -- (-1.5, 0.3) 
        .. controls (-1.6, 0.4) and (-1.6, 0.5) .. (-1.5, 0.6);

    \draw[thick,  fill= red] 
        (-1.5, -0.6) .. controls (-1.2, -0.8) and (-0.8, -1.2) .. (0, -1.2) 
        .. controls (0.8, -1.2) and (1.2, -0.8) .. (1.5, -0.6)
        .. controls (1.6, -0.5) and (1.6, -0.4) .. (1.5, -0.3)
        -- (-1.5, -0.3) 
        .. controls (-1.6, -0.4) and (-1.6, -0.5) .. (-1.5, -0.6);

    \fill[gray!50] (-1.5, 0.3) -- (1.5, 0.3) -- (1.5, -0.3) -- (-1.5, -0.3) -- cycle;
\end{scope}

    \end{tikzpicture}
    }\label{fig:e}
}
& 

\subfigure [] {
  \resizebox{0.3\textwidth}{!}{%
        \begin{tikzpicture}[scale=0.35, baseline={(current bounding box.south)}]
  \draw[thick] (0, 0) circle (5);
\draw[fill=red, thick] 
    (-2, 1.2) .. controls (-0.8, 2) and (0.8, 2) .. (2, 1.2)  
        .. controls (2.5, 0.5) and (-2.5, 0.5) .. (-2, 1.2);
\draw[fill=red, thick] 
    (-1.5, -1.0) .. controls (-0.8, -2) and (0.8, -2) .. (1.5, -1.0) 
        .. controls (2, -0.5) and (-2, -0.5) .. (-1.5, -1.0);
        \draw[thick, ->] (-4.25, 2.7) -- (-0.95, 1.65);
        \draw[thick, ->] (-0.95, 1.65) -- (0.5, 4.975);
    \end{tikzpicture}
    }\label{fig:f}
}
\end{tabular}

\end{figure}

\,

\noindent\textbf{Broken Rays in the Presence of Obstacles.}
When \(\mathcal{E}\) and \(\mathcal{R}\) are separate components of \(\partial M\), the broken rays reflect off obstacles within the domain. 
This case is models the problem where \(M\) contains internal convex obstacles (see Figure \subref{fig:a}). The article \cite[Theorem 3.1]{Eskin:2004} establishes injectivity for the broken ray transform in the presence of several convex obstacles with corner points in \(\mathbb{R}^2\), playing a crucial role in inverse boundary value problems for Schrödinger equations with electromagnetic potentials. However, it is not easy to understand the geometric meaning of all the assumptions posed. 

On Riemannian surfaces with a strictly convex obstacle, injectivity was demonstrated using energy estimates under the assumption of negative curvature \cite{Ilmavirta:Salo:2016}. These results have been generalized to higher dimensions and symmetric tensor fields of any order in \cite{Ilmavirta:Paternain:2022}. In \cite{Jathar:2024}, injectivity up to gauge is established for the broken twisted ray transform involving sums of functions and 1-forms on surfaces with strictly twisted convex obstacles, under nonpositive curvature assumptions. While significant progress has been made, several open problems remain interesting directions to continue research in the future.
\begin{enumerate}[label=(\roman*)]
    \item The question of injectivity in a domain \(\Omega \subset \mathbb{R}^2\) containing two convex smooth obstacles, particularly when one obstacle is strictly convex, remains unresolved \cite[p. 21]{ilmavirta2014thesis} (see Figure \subref{fig:f}). Interestingly, a counterexample exists in the cases where both obstacles are convex, indicating non-injectivity (see Figure \subref{fig:e}). In the usual analysis, multiple tangential reflections as well as trapped broken rays introduce several technical challenges.
    \item Establishing injectivity without curvature assumptions for simple manifolds remains an open problem, despite known results under nonnegative curvature assumptions \cite{Ilmavirta:Salo:2016,Ilmavirta:Paternain:2022} and for the usual geodesic ray transform. Results to this direction could be perhaps obtained by analyzing broken index forms in detail.
    \item Stability results for the broken ray transform in the presence of obstacles do not exist at the moment. These results would be quite interesting additions to the literature. The main difficulties lie in certain regularity questions which need to be resolved.
\end{enumerate}

 \section{Reflection with Vertex}
For a comprehensive discussion on the broken ray transform involving vertices, see \cite[Part 2]{MR4652571}. The V-line transform is defined for piecewise trajectories consisting of two rays emanating from a single vertex point \(x \in \mathbb{R}^n\) (see Figure \subref{fig:g}). Given directions \(v_1, v_2 \in \mathbb{S}^{n-1}\), the V-line transform of a function \(f: \mathbb{R}^n \to \mathbb{R}\), denoted by \(R_{\rm v}f\), is expressed as
\[
R_{\rm v}f(x, v_1, v_2) = c_1 \int_0^{\infty} f(x + t v_1) \, dt + c_2 \int_0^{\infty} f(x + t v_2) \, dt
\]
for \(f \in \mathscr{S}(\mathbb{R}^n)\). We refer to recent works on the V-line transform \cite{Ambartsoumian:2020,Ambartsoumian:2024,Ambartsoumian:Zamindar:2024}.
  Similarly, the cone transform is defined for a function \(f: \mathbb{R}^n \to \mathbb{R}\) over a cone \(C(x, \nu, \beta)\), where \(x\) represents the cone's vertex, \(\nu \in \mathbb{S}^{n-1}\) is the axis of symmetry, and \(\beta \in (0, \pi/2)\) is the half-opening angle (see Figure \subref{fig:h}). The cone transform \(R_{\rm c}f\) is given by
\[
R_{\rm c}f(x, \nu, \beta) = \int_{C(x, \nu, \beta)} f \, ds
\]
for \(f \in \mathscr{S}(\mathbb{R}^n)\). Lastly, the star transform is defined for a set of rays emanating from a vertex \(x \in \mathbb{R}^n\) along directions \(v_1, \ldots, v_m\), forming a ``star" shape (see Figure \subref{fig:i}). The star transform of \(f\), denoted by \(R_{\rm s}f\), is
\[
R_{\rm s}f(x, v_1, \ldots, v_k) = \sum_{i=1}^k \int_0^{\infty} f(x + t v_i) \, dt.
\]

\section{Nonlinear Transforms Related to Broken Rays}

\noindent\textbf{Rigidity Results for the Broken Ray Transform.}
Consider a compact manifold \((M, g)\) of dimension \(n \geq 3\) with boundary \(\partial M\). The {broken scattering relation} is defined as the set of triples \[((x, \xi), (y, \zeta), \tau(\gamma)),\] where \((x, \xi)\) represents the influx data (entry point and direction), \((y, \zeta)\) represents the outflux data (exit point and direction), and \(\tau(\gamma)\) is the travel time of the broken ray \(\gamma\) from its entry to exit. In \cite{Kurylev:2010}, it was shown that the broken scattering relation, combined with boundary data, uniquely determines the isometry type of compact Riemannian manifolds of dimension \(n \geq 3\). This result was extended to reversible Finsler manifolds with convex foliations in \cite{deHoop:2021}, proving that two such manifolds with identical broken scattering relations must be isometric.

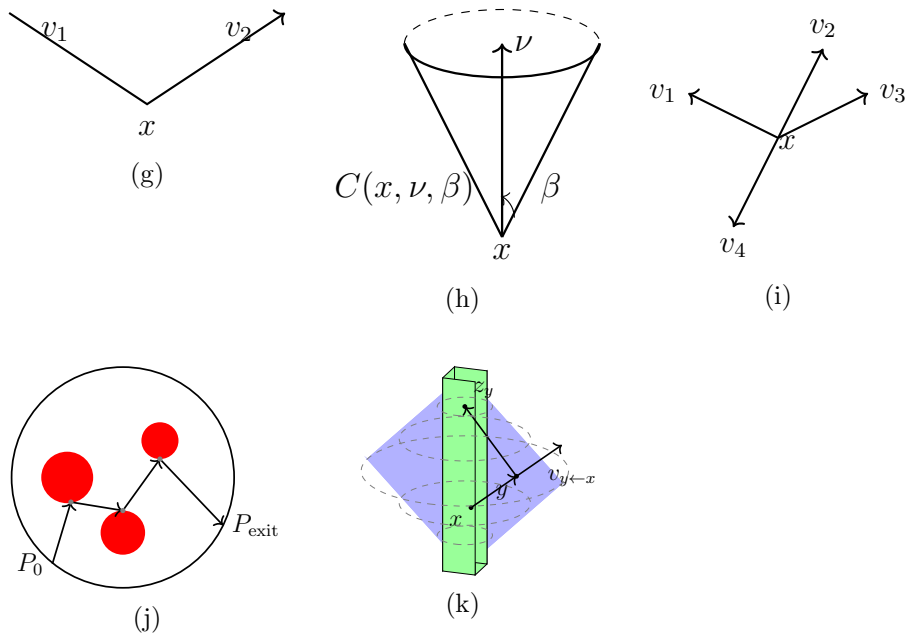
\begin{figure}[htbp]

\begin{tabular}{C{.3\textwidth}C{.3\textwidth}C{.3\textwidth}}
\subfigure [] {
  \resizebox{0.3\textwidth}{!}{%
    \begin{tikzpicture}[scale=0.6, baseline={(current bounding box.south)}]
\coordinate (O) at (0,0);

\draw[thick, ->] (-3,2) -- (O) node[below=2pt] {$x$} -- (3,2);

\node at (-2, 1.6) { $v_1$};
\node at (2, 1.6) {$v_2$};
    \end{tikzpicture}
   }\label{fig:g} 
}& 

\subfigure [] {
  \resizebox{0.3\textwidth}{!}{%
   \begin{tikzpicture}[scale=0.6, baseline={(current bounding box.south)}]
\coordinate (X) at (0,0);

\draw[thick, ->] (X) -- (0,4) node[right] {$\nu$};

\draw[thick] (X) -- (-2,4); 
\draw[thick] (X) -- (2,4);

\draw[thick] (-2,4) arc[start angle=180, end angle=360, x radius=2, y radius=0.7]; 
\draw[dashed] (-2,4) arc[start angle=180, end angle=0, x radius=2, y radius=0.7]; 
\node at (0, -0.3) { $x$};

\draw[->] (0.25,0.4) arc[start angle=0, end angle=63.4, radius=0.5];
\node at (1, 1) { $\beta$};

\node at (-2, 1) { $C(x,\nu,\beta)$};
    \end{tikzpicture}
    }\label{fig:h}
}& 

\subfigure [] {
  \resizebox{0.3\textwidth}{!}{%
    \begin{tikzpicture}[scale=0.6,baseline={(current bounding box.south)}]
   \coordinate (x) at (0,0);

\draw[thick, ->] (x) -- (-2,1) node[left] {$v_1$}; 
\draw[thick, ->] (x) -- (1,2) node[above] {$v_2$};  
\draw[thick, ->] (x) -- (2,1) node[right] {$v_3$};  
\draw[thick, ->] (x) -- (-1,-2) node[below] {$v_4$}; 

\node at (0.2,-0.1) { $x$};
    \end{tikzpicture}
    }\label{fig:i}
}\\

\subfigure [] {
  \resizebox{0.3\textwidth}{!}{%
    \begin{tikzpicture}[scale=0.6, baseline={(current bounding box.south)}]

  \draw[thick] (0,0) circle (3cm);

  \coordinate (ObstacleA) at (-1.5,0);
  \coordinate (ObstacleB) at (1,1);
  \coordinate (ObstacleC) at (0,-1.5);

  \fill[red] (ObstacleA) circle (0.7cm);  
  \fill[red] (ObstacleB) circle (0.5cm);  
  \fill[red] (ObstacleC) circle (0.6cm);

  \coordinate (P0) at (-1.9, -2.35); 
  \coordinate (I1) at (-1.4, -0.67); 
  \coordinate (I2) at (0, -0.9);    
  \coordinate (I3) at (1, 0.5);   
  \coordinate (Exit) at (2.7, -1.3);

  \draw[->, thick] (P0) -- (I1);
  \draw[->,thick] (I1) -- (I2);
  \draw[->,thick] (I2) -- (I3);
  \draw[->,thick] (I3) -- (Exit);

  \fill[gray] (I1) circle (2pt);
  \fill[gray] (I2) circle (2pt);
  \fill[gray] (I3) circle (2pt);

  \node at (P0) [left] {$P_0$};
  \node at (Exit) [right] {$P_{\text{exit}}$};

\end{tikzpicture}
    }\label{fig:j}
}& 

\subfigure [] {
  \resizebox{0.3\textwidth}{!}{%
    \tdplotsetmaincoords{70}{110}

\begin{tikzpicture}[tdplot_main_coords, scale=0.6]

\def\diamondHeight{6}
\def\diamondRadius{3}
\def\stripWidth{1}
\def\coneRadius{1.5}

\coordinate (O) at (0,0,0);

\coordinate (A) at (0,0,\diamondHeight/2);  
\coordinate (B) at (\diamondRadius,0,0);   
\coordinate (C) at (0,\diamondRadius,0);    
\coordinate (D) at (-\diamondRadius,0,0);   
\coordinate (E) at (0,-\diamondRadius,0);   
\coordinate (F) at (0,0,-\diamondHeight/2);

\filldraw[blue!30, opacity=0.4] (A) -- (B) -- (C) -- cycle; 
\filldraw[blue!30, opacity=0.4] (A) -- (C) -- (D) -- cycle; 
\filldraw[blue!30, opacity=0.4] (A) -- (D) -- (E) -- cycle; 
\filldraw[blue!30, opacity=0.4] (A) -- (E) -- (B) -- cycle; 
\filldraw[blue!30, opacity=0.4] (F) -- (B) -- (C) -- cycle; 
\filldraw[blue!30, opacity=0.4] (F) -- (C) -- (D) -- cycle; 
\filldraw[blue!30, opacity=0.4] (F) -- (D) -- (E) -- cycle; 
\filldraw[blue!30, opacity=0.4] (F) -- (E) -- (B) -- cycle;

\draw[fill=green!40, opacity=0.5]
  (-\stripWidth/2,-\stripWidth/2,-\diamondHeight/2) -- (-\stripWidth/2,-\stripWidth/2,\diamondHeight/2)
  -- (\stripWidth/2,-\stripWidth/2,\diamondHeight/2) -- (\stripWidth/2,-\stripWidth/2,-\diamondHeight/2) -- cycle;
\draw[fill=green!40, opacity=0.5]
  (-\stripWidth/2,\stripWidth/2,-\diamondHeight/2) -- (-\stripWidth/2,\stripWidth/2,\diamondHeight/2)
  -- (\stripWidth/2,\stripWidth/2,\diamondHeight/2) -- (\stripWidth/2,\stripWidth/2,-\diamondHeight/2) -- cycle;
\draw[fill=green!40, opacity=0.5]
  (-\stripWidth/2,-\stripWidth/2,-\diamondHeight/2) -- (-\stripWidth/2,\stripWidth/2,-\diamondHeight/2)
  -- (-\stripWidth/2,\stripWidth/2,\diamondHeight/2) -- (-\stripWidth/2,-\stripWidth/2,\diamondHeight/2) -- cycle;
\draw[fill=green!40, opacity=0.5]
  (\stripWidth/2,-\stripWidth/2,-\diamondHeight/2) -- (\stripWidth/2,\stripWidth/2,-\diamondHeight/2)
  -- (\stripWidth/2,\stripWidth/2,\diamondHeight/2) -- (\stripWidth/2,-\stripWidth/2,\diamondHeight/2) -- cycle;

\coordinate (y) at (2.5,2.5,1);
\coordinate (x) at (0.3,0.3,-1);

\coordinate (z_y) at (0,0,\diamondHeight/3);
\fill (z_y) circle (2pt) node[above right] {$z_y$};

\fill (y) circle (2pt) node[below left] {$y$};

\fill (x) circle (2pt) node[below left] {$x$};

\draw[->, thick] (x) -- (y);

\draw[->, thick] (y) -- (z_y);

\draw[->, thick] (y) -- node[below right] {$v_{y \leftarrow x}$} ($(y) + ($(y)-(x)$)$);

\foreach \z in {-2,-1,0,1,2}
{
  \pgfmathsetmacro{\r}{\diamondRadius*(1 - 2.2*abs(\z)/\diamondHeight)}
  \draw[dashed, black!50] (0,0,\z) circle (\r);
}

\end{tikzpicture}
    }\label{fig:k}
}

\end{tabular}

\end{figure}

\noindent\textbf{Broken Non-Abelian Ray Transform.}
Consider a connection \(A\) on the vector bundle \(E = M \times \mathbb{C}^n\), which induces a covariant derivative \(d_A\) on functions \(f: M \to \mathbb{C}^n\), defined by \(d_A f = df + Af\). For a broken ray \(\gamma \in \Gamma\) with a single reflection point \(t^*\) within the manifold \(M\), let \(P_{\gamma[0, t^*]}^A\) and \(P_{\gamma[t^*, \tau]}^A\) denote the parallel transport maps with respect to the connection \(d_A\) along \(\gamma\) before and after the reflection, respectively. The scattering data associated with \(\gamma\) under the connection \(A\) is then given by
\[
S_{\gamma}^A = P_{\gamma[t^*, \tau]}^A P_{\gamma[0, t^*]}^A.
\]
The broken non-abelian ray transform associates the connection \(A\) with this scattering data over all broken rays \(\gamma \in \Gamma\). The broken non-abelian ray transform appears naturally in inverse problems related to the Yang-Mills-Higgs equations \cite{CLOP2,CLOP1}. 

In \cite[Theorem 1.4]{stamant2024gaussian}, it was shown that two smooth Hermitian connections \(A_1\) and \(A_2\) yield the same broken non-abelian ray transform if and only if they are gauge equivalent. This result plays a significant role in inverse problems, particularly in reconstructing connections from the Dirichlet-to-Neumann map at high fixed frequencies. The study in \cite{RefJ6} examined the broken non-abelian ray transform in Minkowski space (see Figure \subref{fig:k}), establishing a stability estimate that accounts for gauge transformations and providing a new proof of injectivity. An application was demonstrated in recovering specific connections, such as light-sink connections, from noisy X-ray transform data using Bayesian inversion techniques.

\section*{Acknowledgements}
We thank Sean Holman and Bill Lionheart, the organizers of the minisymposium M13: Applications of Rich Tomography at the IPMS 2024, Malta, for giving us an opportunity to present our work. The work of S.R.J. and J.R. was supported by the Research Council of Finland through the Flagship of Advanced Mathematics for Sensing, Imaging and Modelling (decision number 359183).

\bibliography{math}

\bibliographystyle{alpha}
\end{document}